\documentclass[11pt]{article}
\usepackage[utf8]{inputenc}
\usepackage[legalpaper,margin=1.0in]{geometry}
\usepackage{authblk}
\usepackage[round]{natbib}
\usepackage{microtype}
\usepackage{graphicx}
\usepackage[caption=false,font=footnotesize]{subfig}
\usepackage{booktabs} 
\usepackage{cite}
\usepackage{apptools}
\usepackage{amssymb,xcolor}
 \usepackage{algorithm} 
\usepackage{algcompatible}
\usepackage{enumitem}
\usepackage[noend]{algpseudocode}
\usepackage{url}
\usepackage{isomath,amsmath,amsthm}
\usepackage{amsfonts,amscd,mathtools}
\usepackage{amssymb}
\usepackage{bm}
\def\D{\bm{D}}
\def\X{\bm{X}}
\def\Pb{\bm{P}}
\def\p{\bm{p}}
\def\real{\mathbb{R}}
\def \alphab {\bm{\alpha}}
\def \betab {\bm{\beta}}
\def \v{\bm{v}}
\def \w{\bm{w}}

\def \alphab {\bm{\alpha}}
\def \betab {\bm{\beta}}
\def \e{\bm{e}}
\def\H{\bm{H}}

\def\A{\bm{A}}
\def\C{\bm{C}}
\def\I{\bm{I}}
\def\J{\bm{J}}
\def\M{\bm{M}}
\def\W{\bm{W}}

\def\a{\bm{a}}
\def\b{\bm{b}}
\def\x{\bm{x}}

\def\Kn{\mathcal{K}_n}
\def\Tn{\mathcal{T}_n}
\def\Adj{\mathcal{A}}

\def\univ{\mathbb{I}}

\def\bg{\bigg}

\def\rmk{\qquad\text}
\def\ones{\bm{1}}

\DeclareMathOperator*{\argmin}{arg\,min}

\newcommand*{\matminus}{%
  \leavevmode
  \hphantom{0}%
  \llap{%
    \settowidth{\dimen0 }{$0$}%
    \resizebox{1.1\dimen0 }{\height}{$-$}%
  }%
}

\begin{document}
\title{A dual basis approach to multidimensional scaling}
\date{\vspace{-4ex}}

\author[1]{Samuel Lichtenberg}
\author[1]{Abiy Tasissa}
\affil[1]{Department of Mathematics, Tufts University, Medford, MA 02155, USA
}

\maketitle
\title{}

\begin{abstract}
Classical multidimensional scaling (CMDS) is a technique that embeds a set of objects in a Euclidean space given their pairwise Euclidean distances. 
The main part of CMDS involves double centering a squared distance matrix and using a truncated eigendecomposition to recover the point coordinates. 
In this paper, motivated by a study in Euclidean distance geometry, we explore a dual basis approach to CMDS. We give an explicit formula for the dual basis vectors and fully characterize the spectrum of an essential matrix in the dual basis framework. We make connections to  a related problem in metric nearness. 
\end{abstract}

\section{Introduction}
The origin of our work is the following inverse problem: given a set of squared Euclidean distances $\{D_{i,j}\}_{1\le i<j\le n}$ among $n$ objects, arranged into a matrix $\D \in \real^{n\times n}$, find points $\p_{1},\p_{2},...,\p_{n}$ in $\real^{r}$ that realize these distances.
The matrix $\Pb \in \real^{n \times r}$ comprises the points arranged as its rows. 
We assume that $n>r$ and the set of points span $\real^{r}$. If all distances are given, the problem can be solved via classical multidimensional scaling (CMDS)  \citep{young1938discussion,torgerson1952multidimensional,torgerson1958theory,gower1966some}. An important result in CMDS states that  
$\D$ is a squared Euclidean matrix if and only if $\X = -\frac{1}{2}\J\D\J$ is positive semidefinite, where $\J=\I - \frac{1}{n}\ones \ones^{\top}$ is the centering matrix \citep{schoenberg1935remarks}. A useful fact is that $\text{rank}(\X)=\text{rank}(\Pb)$. With that, an $r$-truncated eigendecomposition, where $r = \text{rank}(\Pb)=\text{rank}(\X)$, then recovers $\Pb$ from $\X$.

However, if some distances are missing, the problem is studied under the name Euclidean distance geometry (EDG). It is of fundamental importance in many applications \citep{fang2013using,lavor2012recent,ding2010sensor,biswas2006semidefinite}. Building off of matrix completion theory \citep{recht2010guaranteed, candes2009exact}, 
the work in \citep{abiy_exact} analyzes the EDG problem using a dual basis approach. A key part of that work was to represent any zero-centered Gram matrix $\X$ by the expansion 
\begin{equation}\label{eq:x_expansion_dual}
\X = \sum_{\alphab\in\univ} \langle \X\,,\w_{\alphab}\rangle \v_{\alphab},
\end{equation}
where 
$\{ \w_{\alphab} \}$ is a basis of the zero-centered symmetric matrices, $\{ \v_{\alphab} \}$ is its dual basis, and $\univ$ denotes the universal set $\{ (i, j) : 1 \leq i < j \leq n \}$, which has size $L = n(n-1)/2$.
In particular, for $\alphab = (\alpha_1,\alpha_2)$, \citep{abiy_exact} defined
\begin{equation}\label{eq:basis}
\w_{\alphab} = \e_{\alpha_{1},\,\alpha_{1}}+\e_{\alpha_{2},\,\alpha_{2}}-\e_{\alpha_{1},\,\alpha_{2}}-\e_{\alpha_{2},\,\alpha_{1}},
\end{equation}
where $\e_{i,j}$ represents the matrix of zeros except a $1$ at the $(i,j)$-th entry. The matrix 
$\H$ denotes the inner product matrix defined entrywise as $H_{\alphab, \betab} = \langle \w_{\alphab}, \w_{\betab} \rangle $ where $\langle \cdot,\cdot \rangle$ denotes the trace inner product and
\begin{equation}
\v_{\alphab}= \sum_{\betab\in \univ} [\H^{-1}]_{\alphab,\,\betab} \w_{\betab}.
\end{equation}
By construction, $\{\v_{\alphab}\}$ is a dual basis of $\{\w_{\alphab}\}$ satisfying $\langle \v_{\alphab}\,,\w_{\betab} \rangle =
\delta_{\alphab}^{\betab}$ where $\delta_{\alphab}^{\betab}=1$ if $\alphab =\betab$ and zero otherwise. Note that $\alphab = \betab$ holds if $\alpha_1=\beta_1$ and $\alpha_2=\beta_2$. The matrix $\H^{-1}$ is an inner product matrix of the dual basis. Specifically, $[\H^{-1}]_{\alphab,\,\betab}=\langle \v_{\alphab}\,,\v_{\betab} \rangle$. Figure \ref{fig:basis-objects} gives some concrete examples of the core objects that occur in the dual basis approach of \citep{abiy_exact}, for the case $n=4$. Despite the centrality of these objects, the spectrum of $\H$ was not completely characterized, and the only known form for $\v_{\alphab}$ (that did not involve $\H^{-1}$) had multiple cases. 
\begin{figure}
    \centering
    \begin{align*}
    &\w_{1,2} = \begin{bmatrix*}[r]
    1 & \matminus 1 & 0 & 0 \\ 
    \matminus 1 & 1 & 0 & 0 \\
    0 & 0 & 0 & 0 \\
    0 & 0 & 0 & 0
    \end{bmatrix*}
    \hspace{24pt}
    && \v_{1,2} = \frac{1}{16}\begin{bmatrix*}[r]
        3 & \matminus5 & 1 & 1 \\
        \matminus5 & 3 & 1 & 1 \\
        1 & 1 & \matminus 1 & \matminus 1 \\
        1 & 1 & \matminus 1 & \matminus 1
    \end{bmatrix*} \\[12pt]
    &\H = \begin{bmatrix*}[r]
        4 & 1 & 1 & 1 & 1 & 0 \\
        1 & 4 & 1 & 1 & 0 & 1 \\
        1 & 1 & 4 & 0 & 1 & 1 \\
        1 & 1 & 0 & 4 & 1 & 1 \\
        1 & 0 & 1 & 1 & 4 & 1 \\
        0 & 1 & 1 & 1 & 1 & 4
    \end{bmatrix*}
    \hspace{10pt}
    && \H^{-1} = \frac{1}{16}\begin{bmatrix*}[r]
        5 & \matminus 1 & \matminus 1 & \matminus 1 & \matminus 1 & 1 \\
        \matminus 1 & 5 & \matminus 1 & \matminus 1 & 1 & \matminus 1 \\
        \matminus 1 & \matminus 1 & 5 & 1 & \matminus 1 & \matminus 1 \\
        \matminus 1 & \matminus 1 & 1 & 5 & \matminus 1 & \matminus 1 \\
        \matminus 1 & 1 & \matminus 1 & \matminus 1 & 5 & \matminus 1 \\
        1 & \matminus 1 & \matminus 1 & \matminus 1 & \matminus 1 & 5
    \end{bmatrix*}.
    \end{align*} 
    \caption{Dual basis objects $\w_{1,2}, \v_{1,2}, \H$ and  $\H^{-1}$ for $n=4$}
    \label{fig:basis-objects}
\end{figure}

In this work, we look to complete the picture.
We study the exact case, where we have complete distance information in $\D$, and so our work can be considered a dual basis approach to CMDS. However, we characterize and analyze objects that were developed for the EDG problem, and so our results remain applicable to that setting.

\paragraph{Contributions} The contributions of this paper are as follows. First, we characterize all eigenvalues of the matrix $\H$, by relating the structure of $\H$ to a well-studied graph in spectral graph theory \citep{Cameron2003StronglyRG}. The second contribution of this paper is a simple and explicit form for both the dual basis matrices $\{ \v_{\betab} \}$ and their non-zero eigenvalues and corresponding eigenvectors. The third contribution is to express the stability of the map from $\D$ to $\X$ in \eqref{eq:x_expansion_dual} in terms of the entries of $\D$. 
Finally, we connect the matrix $\H$ to a matrix studied in the metric nearness problem; we use our new understanding of $\H$ to give a proof of an empirical observation made in \citep{dhillon_metric_nearness}. Notation used throughout the paper is given in Table \ref{tab:notation}.

\begin{table}[h!]
\centering
\scalebox{0.87}{
    \begin{tabular}{|l|l||l|l|}
        \hline 
         $\x$ & Vector & $\univ$ & Index set $\{ (i,j) : 1 \leq i < j \leq n \}$ \\
              $x_i$ & Vector entry & $\I, \I_n$ & Identity matrix \\ 
         $\X$  & Matrix & $\mathcal{G} $ & Graph \\
         $X_{i,j}, [\X]_{i,j}$ & Matrix entry & $\Kn$ & Complete graph on $n$ vertices \\
         $\X(i,:), \X(:,j)$ & Matrix row, column & $\Adj(\mathcal{G})$ & Adjacency matrix of $\mathcal{G}$ \\ 
         $\alphab, (\alpha_1, \alpha_2), (i,j)$ & Basis index & $\J$ & Centering matrix $\I - \frac{1}{n}\ones \ones^{\top}$\\
         $\w_{\alphab}, \w_{i,j}$ & Basis matrix & $\H$ & Inner product matrix of $\{ \w_{\alphab} \}$ \\
         $\v_{\betab}, \v_{k,l}$ & Dual basis matrix & $\text{vec}(\X)$ & Vectorized matrix  \\
         $\ones$ & All-ones vector & $\mu(\cdot)$ & Mean of vector or matrix \\
         $\e_i$ & Standard basis vector & $\odot$ & Hadamard product\\
         $\e_{\alphab}, \e_{i,j}$ & Standard basis matrix & $\delta_{\alphab}^{\betab}$ & Kronecker delta ($\alphab = \betab$) \\
         $\langle \X \,,\mathbf{Y} \rangle$ & Trace inner product & $\langle \x\,,\mathbf{y}\rangle$ & Dot product       \\ \hline
    \end{tabular}
    }
    \caption{Notation}
    \label{tab:notation}
\end{table}

\section{Spectrum of H}
\label{}

Consider a set of $n$ points $\{\p_i\}_{i=1}^{n}$  in $\real^r$. Let $L = \binom{n}{2}$ denote the number of pairs $(i,j), 1 \leq i < j \leq n$. In \citep{abiy_exact}, the authors defined an inner product matrix $\H \in \real^{L\times L}$ with entries $H_{\alphab, \betab} = \langle \w_{\alphab}, \w_{\betab} \rangle$. They noted that $\H$ is positive definite and symmetric, and, with indices $\alphab = (\alpha_1, \alpha_2), \betab = (\beta_1, \beta_2)$, the entries of $\H$ are
\begin{align*}
    H_{\alphab, \betab} = \begin{cases}
    4 & \text{if $\alpha_1 = \beta_1$ and $\alpha_2 = \beta_2$} \\
    0 & \text{if $\{ \alpha_1, \alpha_2\} \cap \{ \beta_1, \beta_2 \} = \varnothing$} \\
    1 & \text{otherwise.}
    \end{cases}
\end{align*}
For the spectrum of $\H$, the work in \citep{abiy_exact} only gave bounds for the minimum and maximum eigenvalues, despite empirically observing that for any $n$, $\H$ has only three distinct eigenvalues: $2$, $n$, and $2n$. 

In contrast, we are able to completely characterize the spectrum of $\H$. Before we discuss our results, few definitions are in order. A graph is a  mathematical object that consists of a set of  vertices and a set of edges that connect these vertices. A complete graph is a graph in which every pair of its vertices is connected by an edge. The complete graph with $n$ vertices is denoted by $\Kn$. Our result is based on relating the analysis of spectrum of $\H$  
to the well-studied graph $\Tn$, called the triangular graph \citep{Cameron2003StronglyRG}. 
The graph $\Tn$ is formed from the complete graph on $n$ vertices, $\Kn$: 
add a vertex $v_{i,j}$ in $\Tn$ for each edge $e_{i,j}$ in $\Kn$, and connect two vertices $v_{i,j}, v_{k,l} \in \Tn$ with an edge if $e_{ij}$ and $e_{k,l}$ share a common vertex in $\Kn$.

\newtheorem{h_is_triangular_claim}[thm]{Claim}
\begin{h_is_triangular_claim}
    $\H = 4\I_L + \Adj(\Tn)$, where
 $\Adj(\Tn)$ is the adjacency matrix of the triangular graph $\Tn$.
\end{h_is_triangular_claim}
\begin{proof}

We first note that the diagonal entries of $\H$ are all 4 and the off-diagonal entries are either 0 or 1. With that, $\H$ has a natural decomposition as $\H = 4\I_L + \A$, where $\A$ is the adjacency matrix of some graph. We next consider 
consider $\Kn$. In the map from $\Kn$ to $\Tn$, each edge $e_{i, j}$ of $\Kn$ becomes a vertex $v_{i, j}$ of $\Tn$, and two vertices $v_{i, j}, v_{k, l}$ are joined with an edge in $\Tn$ if the distinct edges $e_{i, j}$ and $e_{k, l}$ meet at a vertex in $\Kn$. 
This is equivalent to the condition that $\{ i, j \} \cap \{ k, l \} \neq \varnothing$. 
Now, if we identify the $n$ vertices of $\Kn$ with the points $\{ \p_i \}_{i=1}^{n}$ and the $L$ edges with the basis vectors $\{ \w_{\alphab}\}$ (where $\alphab = (i, j), i < j$), then we see that the vertices of $\Tn$ are the basis vectors $\{ \w
_{\alphab} \}$ and two vertices $\w_{\alphab}, \w_{\betab}$ share an edge in $\Tn$ if and only if $\langle \w_{\alphab}, \w_{\betab} \rangle = 1$. 
Thus $\A$ is exactly the adjacency matrix $\Adj(\Tn)$ of $\Tn$.
\end{proof}

\newtheorem{h_spectrum}[thm]{Corollary}
\begin{h_spectrum}
\label{cor:specH}
The inner product matrix $\H$ has three distinct eigenvalues: $2$, $n$ and $2n$ with multiplicity $L-n$, $n-1$ and $1$ respectively.
\end{h_spectrum}
\begin{proof}
The spectrum of $\Adj(\Tn)$ is known: it has eigenvalues $-2, n-4,$ and $2n-4$, with respective multiplicities $L - n, n - 1$ and 1 \citep{hoffman_1959}. Since we have that $\H = 4\I_L + \Adj(\Tn)$,
the result follows.
\end{proof}

\section{Representation of the dual basis}

The work \citep{abiy_exact} gave an expression for the dual basis vectors $\{ \v_{\alphab} \}$, but it was presented with multiple cases. This motivates our first result of this section: a simple dyadic form for $\v_{\alphab}$.

\newtheorem{form_of_dual}[thm]{Claim}
\begin{form_of_dual}
\label{thm:form_of_dual}
For any $\alphab \in \univ$, $\v_{\alphab}$ has the form
\begin{align*}
    \v_{\alphab} &= -\frac{1}{2}\bigg(\a \b^{\top} + \b \a^{\top}\bigg),
\end{align*}
where $\a = (\mathbf{e_i} - \frac{1}{n}\mathbf{1}) = \J(:,i)$ and $\b = (\mathbf{e_j} - \frac{1}{n}\mathbf{1}) = \J(:,j)$.
\end{form_of_dual}
\begin{proof}
To show this is the dual basis, it suffices to certify biorthogonality. Specifically, for $\alphab = (i,j)$ and $\betab = (k,l)$, with $i < j$ and $k < l$, we must show that $\langle \w_{\alphab}, \v_{\betab} \rangle = \delta_{\alphab}^{\betab}$. Let $\a = (\e_k - \frac{1}{n}\ones)$ and $\b = (\e_l - \frac{1}{n}\ones)$. We have
\begin{align*}
    \langle \w_{\alphab}, \v_{\betab} \rangle &= -\frac{1}{2}\langle \w_{i,j}, \a \b^{\top} + \b \a^{\top} \rangle \\
    &= -\frac{1}{2}\langle \e_{i,i} + \e_{j,j} - \e_{i,j} - \e_{j,i}, \a \b^{\top} + \b \a^{\top} \rangle \\
    &=  -\frac{1}{2}\bg([\a\b^{\top} + \b \a^{\top}]_{i,i} + [\a\b^{\top} + \b \a^{\top}]_{j,j} - 2[\a\b^{\top} + \b \a^{\top}]_{i,j} 
    \bg) \\
    &= -\frac{1}{2} \bg( 2a_ib_i + 2a_jb_j - 2(a_ib_j + b_ia_j)  \bg) \\
    &= a_ib_j + a_jb_i - a_ib_i - a_jb_j \\
    &= \hspace{10pt} \big(\delta_i^k - \frac{1}{n}\big)\big(\delta_j^l - \frac{1}{n}\big) 
    + \big(\delta_j^k - \frac{1}{n}\big)\big(\delta_i^l - \frac{1}{n}\big) \\
    &\hspace{12pt} - \big(\delta_i^k - \frac{1}{n}\big)\big(\delta_i^l - \frac{1}{n}\big) 
    - \big(\delta_j^k - \frac{1}{n}\big)\big(\delta_j^l - \frac{1}{n}\big) \\
    &= \hspace{10pt} \delta_i^k\delta_j^l - \frac{1}{n}\delta_i^k - \frac{1}{n}\delta_j^l + \frac{1}{n^2} 
    + \delta_j^k\delta_i^l - \frac{1}{n}\delta_i^l - \frac{1}{n}\delta_j^k + \frac{1}{n^2} \\
    & \hspace{12pt} -\delta_i^k\delta_i^l + \frac{1}{n}\delta_i^l + \frac{1}{n}\delta_i^k - \frac{1}{n^2} 
    -\delta_j^k\delta_j^l + \frac{1}{n}\delta_j^k + \frac{1}{n}\delta_j^l - \frac{1}{n^2} \\
    &= \delta_i^k\delta_j^l + \delta_j^k\delta_i^l \\
    &= \delta_i^k\delta_j^l \rmk{(since $i < j, k < l$)} \\
    &= \delta_{\alphab}^{\betab},
\end{align*}
which is the desired biorthogonality.
\end{proof}
From this representation, we can obtain the spectrum of each $\v_{\alphab}$.

\newtheorem{eig_of_dual}[thm]{Claim}
\begin{eig_of_dual}
\label{thm:eig_of_dual}
For any $\alphab \in \univ$, $\v_{\alphab}$ is a rank 2 matrix with nonzero eigenvalues $\frac{1}{2}$ and $-\frac{1}{2} + \frac{1}{n}$ and corresponding eigenvectors $\a - \b$ and $\a + \b$, independent of dimension.
\end{eig_of_dual}
\begin{proof}
For now, we ignore the $-\frac{1}{2}$ factor in $\v_{\alphab}$. We observe that $(\a \b^{\top} + \b \a^{\top}) = \C\A$, where $\C$ is a matrix whose columns are $\a$ and $\b$ respectively and $\A$ is a matrix whose rows are $\b^{\top}$ and $\a^{\top}$ respectively. $\A$ has rank $2$, and since $\a$ and $\b$ are always linearly independent, $\text{rank}(\C \A) = \text{rank}(\A) = 2$. To find the nonzero eigenvalues, we guess the eigenvectors to be $\a + \b$ and $\a - \b$. Noting that $\|\a\|^2 = \|\b\|^2 = \frac{n - 1}{n}$ and $\a^{\top}\b = -\frac{1}{n}$, this gives us eigenvalues of $-\frac{1}{2} + \frac{1}{n}$ and $\frac{1}{2}$, respectively.

\end{proof}

To illustrate the utility of this representation, we immediately use it to prove the following result. We know that if $\X$ is a Gram matrix, then $\langle \X, \w_{\alphab} \rangle$, $\alphab = (i, j)$, is a distance $D_{i,j}$; but if $D_{i,j}$ are distances, is it always the case that $\sum_{\alphab} D_{i,j} \v_{\alphab}$ is a Gram matrix? We show that this is indeed true.

\newtheorem{sum_is_gram}[thm]{Claim}
\begin{sum_is_gram}
\label{thm:sum_is_gram}
Given distances $D_{i,j}$, $
    \X = \sum_{(i,j) \in \univ} D_{i,j} \v_{i, j}$
is a Gram matrix.
\end{sum_is_gram}
\begin{proof}
We have that
\begin{align*}
    \X  &= \sum_{(i,j) \in \univ} D_{i,j} \v_{i,j} \\
    &= -\frac{1}{2}\sum_{(i,j) \in \univ} D_{i,j} \big(\J(:,i) \J(:,j)^{\top} + \J(:,j) \J(:,i)^{\top}\big) \\
    &= -\frac{1}{2} \sum_{i \neq j} D_{i,j} (\J(:,i) \J(:,j)^{\top}) \rmk{(distances are symmetric)} \\
    &= -\frac{1}{2} \sum_{i, j} D_{i,j} (\J(:,i)\J(:,j)^{\top}) \rmk{(self-distance is 0)}.
\end{align*}
Letting $\D := [D_{i,j}]$, we consider a particular entry $(a,b)$ of $\X$:
\begin{align*}
    X_{a,b} &= -\frac{1}{2} \bg[\sum_{i, j} D_{i,j} (\J(:,i)\J(:,j)^{\top})\bg]_{a,b} \\
    &= -\frac{1}{2} \sum_i \sum_j D_{i,j}J_{a,i}J_{j,b} \\
    &=  [-\frac{1}{2} \J \D \J]_{a,b}.
\end{align*}
We see $\X = -\frac{1}{2}\J\D\J$, and so $\X$ is a Gram matrix. 

\end{proof}

\section{Stability of Gram matrix under additive noise}

In many applications, the observed distances contain measurement errors. 
Here, we study a simple model where a true squared distance matrix $\D$ is corrupted with additive noise $\tilde{\D}$. We prove the following result on the stability of the underlying Gram matrix.

\newtheorem{noise_bound_claim}[thm]{Claim}

\begin{noise_bound_claim}
Let $\D_{\text{noisy}} = \D+\tilde{\D}$, $\X$ denote the true Gram matrix, and $\bar{\X}$ be the Gram matrix corresponding to $\D_{\text{noisy}}$.
Then $\|\bar{\X} - \X\|_{\infty} < 4\|\tilde{\D}\|_{\infty}$, where $\|\cdot\|_{\infty}$ gives the largest entry of a matrix in absolute value.
\end{noise_bound_claim}
\begin{proof}
Using the dual basis expansion, we have that 
\begin{align*}
    \|\bar{\X} - \X\|_{\infty} &= \max_{a,b} \bg| \sum_{(i,j) \in \univ} 
    [\v_{i,j}]_{a,b}\tilde{D}_{i,j}\bg| \\
    &= \max_{a,b} \bg| -\frac{1}{2} 
    \sum_{(i,j) \in \univ}
    (J_{a,i}J_{j,b}+ J_{a,j}J_{i,b})
    \tilde{D}_{i,j}\bg| \\
    &= \frac{1}{2}\bg(\max_{a,b} \bg| \sum_{(i,j) \in \univ} 
    (J_{a,i}J_{j,b}+ J_{a,j}J_{i,b})
    \tilde{D}_{i,j}\bg|\bg) \\
    &\leq \frac{1}{2}\bg(\max_{a,b} \sum_{(i,j) \in \univ} 
    |J_{a,i}J_{j,b}|
    |\tilde{D}_{i,j}|
    + \max_{a,b}
    \sum_{(i,j) \in \univ} 
    |J_{a,j}J_{i,b}|
    |\tilde{D}_{i,j}|
    \bg) \\
    &\leq \frac{1}{2}
    \|\tilde{\D}\|_{\infty}
    \bg(
    \max_{a,b}
    \sum_{(i,j) \in \univ}
    |J_{a,i}J_{j,b}| 
    +
    \max_{a,b}
    \sum_{(i,j) \in \univ}
    |J_{a,j}J_{i,b}|
    \bg)
    \\
    &= 
    \|\tilde{\D}\|_{\infty}
    \bg(
    \max_{a,b}
    \sum_{(i,j) \in \univ}
    |J_{a,i}J_{j,b}|    
    \bg).
    \rmk{(symmetry of $J$)}
\end{align*}
The expression $|J_{a,i}J_{j,b}|$ can take on three possible values: $((n-1)/n)^2$, $(n-1)/n^2$, and $1/n^2$.
There is at most one term with value $((n-1)/n)^2$, when $a < b$ and $i = a, j = b$. 
Any term with exactly one of $i=a$ or $j=b$ will have $|J_{a,i}J_{j,b}| = (n-1)/n^2$, and there is a strict upper bound of $2n$ on the number of such terms that can occur. 
Finally, there are $\max\{L-2n-1, 0\}$ terms remaining in the sum over $(i,j) \in \univ$, which must have value $1/n^2$. Thus we get the bound
\begin{align*}
    \|\tilde{\D}\|_{\infty}
    \bg(
    \max_{a,b}
    \sum_{(i,j) \in \univ}
    |J_{a,i}J_{j,b}|    
    \bg)  
    &<
    \|\tilde{\D}\|_{\infty}   
    \bg[
    \bg(
    \frac{n-1}{n}
    \bg)
    ^2
    +
    2n
    \bg(
    \frac{n-1}{n^2}
    \bg) \\
    &\hspace{48pt}
    +
    \max\{L-2n-1, 0\}
    \frac{1}{n^2}
    \bg]
    \\
    &< 4\|\tilde{\D}\|_{\infty}.
\end{align*}
\end{proof}

\section{The Metric Nearness Problem}

In \citep{dhillon_metric_nearness}, \citep{sra_metric_nearness}, and \citep{brickell_metric_nearness}, the authors studied a problem related to but distinct from EDG, called the \emph{metric nearness problem}. In this section, we explain and prove an empirical observation of \citep{dhillon_metric_nearness}.

Some background is first required. Let $\mathcal{D}_{n}$ denote the set of ``dissimilarity matrices'': non-negative, symmetric matrices with zero diagonal, and $\mathcal{M}_n$ be the set of distance matrices: matrices in $\mathcal{D}_{n}$ whose entries satisfy the triangle inequality. Given an input matrix $\D \in \mathcal{D}_{n}$ and a weight matrix $\W \in \mathcal{D}_{n}$, the metric nearness problem seeks to find a nearest valid distance matrix $\M$, in the sense that
\begin{align*}
    \M = \underset{X \in \mathcal{M}_{n}}{\argmin}\,\, \| \W \odot (\X - \D) \|,
\end{align*}
for some norm, where $\odot$ is the Hadamard (elementwise) product.

As part of their approach to this problem, the authors in \citep{dhillon_metric_nearness} defined a tall matrix $\A$ that encodes the constraints enforced by the triangle inequality. 
Specifically, $\A$ has $3 \binom{n}{2}$ rows and $\binom{n}{2}$ columns: each column represents a distance between two points, and each row represents a triangle inequality constraint among a set of three distances: $D_{i,k} - D_{i,j} - D_{j,k} \leq 0$. 
The entries in $\A$ are the sign of each distance in a particular inequality. 

It is easier to see with an example, here among four points $\p_1, \p_2, \p_3, \p_4$. 
We label the distance $D_{i,j}$ as $(\p_i, \p_j)$ and the inequality
$D_{i,k} - D_{i,j} - D_{j,k} \leq 0$ as
$(\p_i, \p_k, \p_j)$; the first two points in the latter list refer to the index of $\D$ with positive sign in the inequality.  
Our $18 \times 6$ matrix $\A$ (abbreviated for brevity) is given in Table \ref{tab:nearness-example}.
\begin{table}
    \centering
    \begin{tabular}{ |c|c|c|c|c|c|c| } 
 \hline
 & $(\p_1, \p_2)$ & $(\p_1, \p_3)$ & $(\p_1, \p_4)$ & $(\p_2,\p_3)$ & $(\p_2,\p_4)$ & $(\p_3,\p_4)$ \\ 
 \hline
 $(\p_1,\p_2,\p_3)$ & 1 & \matminus 1 & 0 & \matminus 1 & 0 & 0\\
 $(\p_1,\p_3,\p_2)$ & \matminus 1 & 1 & 0 & \matminus 1 & 0 & 0 \\
 $(\p_2,\p_3,\p_1)$ & \matminus 1 & \matminus 1 & 0 & 1 & 0 & 0 \\
 $\vdots$ & $\vdots$ & $\vdots$ & $\vdots$ & $\vdots$ & $\vdots$ & $\vdots$ \\
 $(\p_2, \p_3, \p_4)$ & 0 & 0 & 0 & 1 & \matminus 1 & \matminus 1 \\
 \hline
\end{tabular}
    \caption{Triangle inequality constraint matrix $\A$ (four points).}
    \label{tab:nearness-example}
\end{table}

The authors in \citep{dhillon_metric_nearness} made the empirical observation that $\A$ has three distinct singular values: $\sqrt{3n-4}, \sqrt{2n-2}, \sqrt{n-2}$ with multiplicities $n(n-3)/2$, $(n\matminus 1)$, and $1$. However, they did not provide a proof of this fact except for the largest singular value. This motivates our following result.

\newtheorem{metric_nearness_claim}[thm]{Claim} 
\begin{metric_nearness_claim} 
$\A^{\top}\A$ and $-\H$ are identical up to diagonal scaling.
\end{metric_nearness_claim}
\begin{proof}  
Let $L = \binom{n}{2}$, and let $\mathbb{T}$ refer to the set of all tuples $(p_i, p_j, p_k)$ for all $i,j,k \in [n]$ (i.e., the row indices of $\A$). We are interested in the entries of the $L \times L$ matrix $\A^{\top}\A$. 

The diagonal entries are easier to compute. Given any $\alpha = (p_i, p_j), i < j$, we have
\begin{align*}
    [\A^{\top}\A]_{\alpha,\alpha} &= \sum_{t \in \mathbb{T}} A_{t, \alpha}^2 = 3(n-2),
\end{align*}
as this is just the number of nonzero entries in column $\alpha$, which was already given by \citep{dhillon_metric_nearness}. The off-diagonal entries require more analysis. For some $\alpha = (p_i, p_j), \beta = (p_k, p_l)$, we are interested in
\begin{align*}
    [\A^{\top}\A]_{\alpha,\beta} &= \sum_{t \in \mathbb{T}} A_{t,\alpha}A_{t,\beta}.
\end{align*}
If none of the constituent points of $\alpha$ and $\beta$ are the same, then the entry $[\A^{\top}\A]_{\alpha, \beta}$ is 0. 
If the first point is shared, such that 
$\alpha = (p_i, p_j), \beta = (p_i, p_k)$, 
then the rows that will have nonzero 
entries are $t_1 = (p_j, p_k, p_i), t_2 = 
(p_i, p_j, p_k), t_3 = (p_i, p_k, p_j)$. We 
have $A_{t_1,\alpha} = \matminus 1, A_{t_1,\beta} = \matminus 1$; $A_{t_2,\alpha} = 1, A_{t_2,\beta} = 
\matminus 1$; $A_{t_3, \alpha} = \matminus
1, A_{t_3, \beta} = 1$. 
So the sum is \matminus 1. And, in fact, this is true whether it is the first point
or the second point of $\alpha$ and $\beta$ that are shared. So in this case, 
$[\A^{\top}\A]_{\alpha, \beta} = \matminus 1$. 

From this, we can see that $\A^{\top}\A$  captures the same structure as the inner product matrix $\H$ that arises from the dual basis approach.
Specifically, we have that
$\A^{\top}\A = (3n-2)\I_L - \H$.
\end{proof}
With this result and Corollary \ref{cor:specH}, we see that the eigenvalues of $\A^{\top}\A$ are $3n-4, 2n-2, n-2$ with the multiplicities $L - n = n(n-3)/2$, $n-2$, and $1$, which, after taking square roots, exactly matches the result of \citep{dhillon_metric_nearness}.

\section{Acknowledgements}
This work is supported by NSF DMS 2208392.

\section{Conclusion}

In this paper, motivated by a previous work in distance geometry, we studied a dual basis framework for classical multidimensional scaling (CMDS). 
We characterized the spectrum of the inner product matrix $\H$, gave a simple form for the dual basis vectors $\v_{\alpha}$, analyzed the stability of an important map, and explained an empirical observation of a related work in metric nearness. 
Our work considered only the exact case, so an important direction for future work is to develop more theory for the sampled regime. 
In addition, we surmise that further theoretical analysis of the spectrum of $\X$, the Gram matrix, can be done using tools from the theory of matrix pencils  \citep{ikramov1993matrix}. 

\bibliographystyle{IEEEtranN}
\bibliography{refs}
\end{document}